\documentclass[11pt]{article}

\usepackage{amsfonts}
\usepackage{amssymb,amsmath,amsthm}
\usepackage{latexsym}
\usepackage{fullpage}
\usepackage{hyperref}
\usepackage{graphicx}
\usepackage{tkz-euclide,subfigure}
\usepackage{tikz}
\usepackage[utf8]{inputenc}
\usepackage{pgfplots}
\usetikzlibrary{shapes.misc}
\usepackage{ amssymb }

\theoremstyle{definition}

\newcounter{tenumerate}

\renewcommand{\epsilon}{\varepsilon}

\newcommand{\remove}[1]{}

\renewcommand{\leq}{\leqslant}
\renewcommand{\geq}{\geqslant}

\def\XXint#1#2#3{{\setbox0=\hbox{$#1{#2#3}{\int}$}
\vcenter{\hbox{$#2#3$}}\kern-.5\wd0}}

\title{\textbf{A class of quadratic reflected BSDEs with singular coefficients}}
\begin{document}
\author{{Shiqiu Zheng\thanks{Corresponding author}\ , \ Lidong Zhang, \ Xiangbo Meng}\\
  \small College of Sciences, Tianjin University of Science and Technology, Tianjin 300457, China
  \\
\small \emph{E-mail: shiqiu@tust.edu.cn}}
\date{}
\maketitle
\begin{abstract}
In this paper, we study the existence and uniqueness of the solution to a reflected backward stochastic differential equation (RBSDE) with the generator $g(t,y,z)=G_f^F(t,y,z)+f(y)|z|^2$, where $f(y)$ is a locally integrable function defined on an open interval $D$, and $G_f^F(t,y,z)$ is induced by $f$ and a Lipschitz continuous function $F$. Both the solution $Y_t$ and the obstacle $L_t$ of this RBSDE take values in $D$. As applications, we provide a probabilistic interpretation of an obstacle problem for a quadratic PDE with a singular term, whose solution takes values in $D$, and study an optimal stopping problem for the payoff of American options under general utilities.
\end{abstract}
\section{Introduction}
El Karoui et al. \cite{E} introduced the notion of reflected backward stochastic differential equations (RBSDEs) with Lipschitz continuous generators. Several works have been done to study RBSDEs whose generator $g(t,y,z)$ is continuous in $(y,z)$ and has a quadratic growth in the variable $z$ (quadratic RBSDEs). The existence and uniqueness of quadratic RBSDEs have been investigated by Kobylanski et al. \cite{K2} and Xu \cite{Xu} for bounded terminal variables and obstacles, and by Bayraktar and Yao \cite{BY} for exponentially integrable terminal variables and obstacles. When the generator is differentiable, the existence and uniqueness of quadratic RBSDEs driven by a continuous martingale were obtained by Lionnet \cite{Li} for bounded terminal variables and upper-bounded obstacles.

In the literature, the generator $g(\cdot,y,z)$ of the quadratic RBSDE is usually assumed to be continuous on $\textbf{R}\times \textbf{R}^d$. This means that such an RBSDE may not be applicable to some problems described by systems involving singular coefficients, such as the obstacle problem involving the following semilinear PDE with a singular term:
$$\partial_tv+\frac{1}{2}\Delta v+h(v)|\nabla v|^2=0,\eqno(1.1)$$
where $h(\cdot)$ is a continuous function defined on $(0,\infty)$.
The existence of a bounded Sobolev solution of an obstacle problem related to (1.1) was studied by Arcoya et al. \cite{A} in the elliptic case. However, to the best of our knowledge, the viscosity solution of the obstacle problem for (1.1) has not been obtained in the literature. Motivated by this, we will study the existence and uniqueness of the RBSDE with generator: $$g(t,y,z)=G^F(t,y,z)+f(y)|z|^2,\ \ \textrm{where}\ \  G^F(t,y,z):=\frac{F(t,u_f(y),u_f'(y)z)}{u_f'(y)},\eqno(1.2)$$
where $f(y)$ is a locally integrable function defined on an open interval $D$,
$$u_f(x):=\int^x_\alpha\exp\left(2\int^y_\alpha f(z)dz\right)dy,\ \ x\in D,$$
and $F$ is a Lipschitz continuous function. Both the solution $Y_t$ and obstacle $L_t$ of such RBSDE take values in $D$. Some typical examples of $u_f$ and $G^F$ are given in Example 2.1. When $f(y)=\frac{1}{y}$, some related BSDEs were studied by Bahlali and Tangpi \cite{BT}. Our study is inspired by Bahlali et al. \cite{B17} and can be seen as an extension of a study by Zheng et al. \cite{Zheng1}. Our proof of is based on the transformation $u_f(y)$, an It\^{o}-Krylov formula, and an existence and uniqueness theorem for an RBSDE whose solution and obstacle take values in an open interval (see Proposition 3.2). The transformation $u_f(y)$ and the It\^{o}-Krylov formula are used to remove the quadratic term $f(y)|z|^2$, where $f(y)$ may be discontinuous. The transformation $u_f(y)$ has been applied to  stochastic differential utility by Duffie and Epstein \cite{DE} and to quadratic BSDEs by \cite{B17} (see also \cite{B19, BT, Zheng1}). In our situation, it is crucial to guarantee that the solution $Y$ takes values in $D$. To this end, we establish an existence and uniqueness result for an RBSDE with a Lipschitz continuous generator, whose solution $Y_t$ and obstacle $L_t$ take values in an open interval. Some typical examples of such RBSDEs are also provided; in particular, we present an example showing that the quadratic BSDE with a bounded terminal variable may have no bounded solution, even when $D=\textbf{R}$, $f(y)$ is a constant, and the corresponding RBSDE has a bounded solution (see Example 3.9).

As an application, we study an obstacle problem for a PDE that is more general than (1.1). By using the comparison theorem established in this study and the relation between the RBSDE with generator $G_f^F+f(y)|z|^2$ and the RBSDE with generator $F$, we provide a viscosity solution to this obstacle problem.

Interestingly, we observe that the transformation $u_f(y)$ can be viewed as a utility function induced by $f(y)$, which includes exponential utility, power utility, and logarithmic utility as special cases. Moreover, $-2f(y)$ represents the Arrow--Pratt risk-aversion coefficient of the utility function $u_f(y)$. Based on this observation, we apply the RBSDE with generator $G_f^F+f(y)|z|^2$ to study an optimal stopping problem for the payoff of American options under the utility function $u_f(y)$.

The remained of the paper is organized as follows. Section 2 presents some assumptions. Section 3 investigates the RBSDE with generator (1.2). Section 4 presents two applications. The Appendix contains auxiliary results.
%%%%%%%%%%%%%%%%%%%%%%%%%%%%%%%%%%%%%%%%%%%%%%%%%%%%%%%%%%%%%%%%%%%%%%%%%%%%%%%%%%%%%%%%%%%%%%%%%%%%%%%%%%%%%%%%%%%%%%%%%%%%%%%%%%%%%%%%%%%%%%%%%%%%%%%%%%%%%%
\section{Preliminaries}
Let $(\Omega ,\cal{F},\mathit{P})$ be a complete probability space. Let ${{(B_t)}_{t\geq
0}}$ be a $d$-dimensional standard Brownian motion defined on this probability space, and let $({\cal{F}}_t)_{t\geq 0}$ be
the natural filtration generated by ${{(B_t)}_{t\geq 0}}$, augmented
by the $\mathit{P}$-null sets of ${\cal{F}}$. Let ${\cal{P}}$ be the progressively measurable sigma-field on $[0,T]\times\Omega.$ For $z\in {\textbf{R}}^d$, let $|z|$ denote the
Euclidean norm. Let $T>0$ and $p>1$ be given real numbers, and let ${\cal{T}}_{t,T}$ be the set of all stopping times $\tau$ satisfying $t\leq \tau\leq T$. Throughout, we assume that $D\subset \textbf{R}$ is an open interval. We define the following spaces:

$L_{1,loc}(D)=\{f:\ f:D\rightarrow \textbf{R},$ is measurable and locally integrable$\}$;

$L_D({\cal{F}}_T)=\{\xi:\ {\cal{F}}_T$-measurable random variable whose range is included in $D\}$;

$L^r_D({\mathcal {F}}_T)=\{\xi\in L_D({\cal{F}}_T)$: ${{E}}\left[|\xi|^r\right]<\infty\},\ r\geq1; $

$L^\infty_D({\mathcal {F}}_T)=\{\xi\in L_D({\cal{F}}_T)$: the range of $\xi$ is included in a closed subset of $D\}$;

${\mathcal{C}}_D=\{(\psi_t)_{t\in[0,T]}:$ continuous and $({\cal{F}}_t)$-adapted process whose range is included in $D\}$;

${\mathcal{S}}^r_D=\{(\psi_t)_{t\in[0,T]}\in{\mathcal{C}}_D$: ${{E}}\left
[{\mathrm{sup}}_{0\leq t\leq T} |\psi _t|^r\right] <\infty \},\ r\geq1;$

${\mathcal{S}}^\infty_D=\{(\psi_t)_{t\in[0,T]}\in{\mathcal{C}}_D$: the range of $(\psi_t)_{t\in[0,T]}$ is included in a closed subset of $D\}$;

${\mathcal{A}}=\{(\psi_t)_{t\in[0,T]}:$ increasing, continuous and  $({\cal{F}}_t)$-adapted
$\textbf{R}$-valued process with $\psi_0=0$$\}$;

${\mathcal{A}}^r=\{(\psi_t)_{t\in[0,T]}\in{\mathcal{A}}$: ${{E}}
[|\psi _T|^r] <\infty\},\ r>0;$

${H}^r=\{(\psi_t)_{t\in[0,T]}:$  $\textbf{R}^d$-valued, $({\cal{F}}_t)$-progressively measurable and $\int_0^T|\psi_t|^rdt
<\infty \},\ r\geq1;$

${\cal{H}}^r=\{(\psi_t)_{t\in[0,T]}\in{H}^2$: ${{E}}[(\int_0^T|\psi_t|^2dt)^{\frac{r}{2}}]<\infty \},\ r\geq1;$

$W^2_{1,loc}(D)=\{f\in L_{1,loc}(D):$ its generalized derivation $f'$ and $f''$ both belong to $L_{1,loc}(D)\}$.\\
For convenience, when the range of a random variable (or a process) is $\textbf{R}$ or is clear from context, we use the simplified notations: $L({\mathcal {F}}_T), L^r({\mathcal {F}}_T), {\mathcal{C}}$, and ${\mathcal{S}}^r.$ Note that in this paper, all the equalities and inequalities for random variables are understood to hold in the almost sure sense.

Let $f\in L_{1,loc}(D)$ be given. For $\alpha\in D,$ we define the transformation:
$$u^\alpha_f(x):=\int^x_\alpha\exp\left(2\int^y_\alpha f(z)dz\right)dy,\ \ x\in D.$$
We assume that $\alpha\in D$ is given and denote $u^\alpha_f(x)$ by $u_f(x)$. Let $V:=\{y:y=u_f(x),x\in D\}.$ Some properties of $u_f(x)$ are provided in Lemma A.1 in the Appendix. Since $u_f(x)$ is continuous, strictly increasing and $u_f(\alpha)=0$, it follows that $V$ is an open interval such that $0\in V.$ In particular, when $f\geq0$ and $D=(a,\infty)$, since $u_f(x)\geq u_0(x)=x-\alpha$ (by Lemma A.1(v)), we have $V=(b,\infty)$ for some $b\geq-\infty.$

Let $\delta\geq0$, $\gamma\geq0$, and $\kappa\geq0$ be given constants. Define a function
$${F}\left(\omega,t,y,z\right): \Omega\times[0,T]\times\mathbf{R}\times\mathbf{R}^{\mathit{d}}\longmapsto \mathbf{R},$$ such that $F$ is measurable with respect to ${\cal{P}}\otimes{\cal{B}}({\mathbf{R}})\otimes{\cal{B}}({\mathbf{R}}^d)$. We always assume that $F$ satisfies\\\\
\textbf{Assumption (A)}. For each $t\in[0,T]$ and $(y,z), (\tilde{y},\tilde{z})\in \textbf{R}{\mathbf{\times R}}^{\mathit{d}},$ we have
  $$|F(t,0,0)|\leq\delta\ \ \textmd{and} \ \ |{F}(t,y,z)-{F}(t,\tilde{y},\tilde{z})|\leq \gamma|y-\tilde{y}|+\kappa|z-\tilde{z}|.$$

We define the function $G_f^F(t,y,z)$ as follows
$$G_f^F(t,y,z):=\frac{F(t,u_f(y),u_f'(y)z)}{u_f'(y)},\ \ \ (t,y,z)\in[0,T]\times D{\mathbf{\times R}}^{\mathit{d}}.$$
Since $F$ satisfies Assumption (A), it follows that $G_f^F$ is continuous in $(y,z)$ and when $\delta=\gamma=0$, $|G_f^F(t,y,z)|\leq\kappa|z|$. Moreover, if $f\geq0,$ then $u_f'(y)$ is nondecreasing, and by the definitions of $G_f^F$ and $u_f$, we can check that for each $c\in D$ and $(t,y,z)\in[0,T]\times([c,\infty)\cap D){\mathbf{\times R}}^{\mathit{d}},$
\begin{eqnarray*}
\ \ \ \ \  \ \ \ \ \ \ \ \  \ \ \ |G_f^F(t,y,z)|&\leq&\frac{\delta}{u'_f(y)}+\frac{\gamma|\int_\alpha^y u'_f(x)dx|}{u'_f(y)}+\kappa|z|\\
&\leq&\frac{\delta}{u'_f(c)}+\frac{\gamma|\int^c_\alpha u'_f(x)dx|}{u'_f(c)}+\gamma\left|\int_c^y\frac{u'_f(x)}{u'_f(y)}dx\right|+\kappa|z|\\
&\leq&\frac{\delta}{u'_f(c)}+\frac{\gamma|u_f(c)|}{u'_f(c)}+\gamma|c|+\gamma|y|+\kappa|z|.\ \ \ \ \  \ \ \ \ \ \ \ \ \ \  \ \ \ \ \ \ \ \  \ \ \ \ \ \ \ (2.1)
\end{eqnarray*}

Let $\beta\in\textbf{R}$, and $\delta_1, \gamma_1$, $\kappa_1\in\textbf{R}$ such that $|\delta_1|\leq\delta, |\gamma_1|\leq\gamma$, and $|\kappa_1|\leq\kappa.$ We list some typical examples of the functions $u_f$ and $G_f^F$.\\\\
\textbf{Example 2.1} \emph{(i) Let $F(t,a,b)=\kappa_1b$ (or $F(t,a,b)=\kappa_1|b|$). Then, for each $f\in L_{1,loc}(D)$,
$$G_f^F(t,y,z)=\kappa_1z\ (or \ G_f^F(t,y,z)=\kappa_1|z|).$$}
\ \ \ \ \emph{(ii) Let} $D=\textbf{R}$ \emph{and $f(y)=0$. Then} $$u_f(y)=y-\alpha\ \ and \ \ G_f^F(t,y,z)=F(t,y-\alpha,z).$$

\emph{(iii) Let $D=(0,\infty),$ $f(y)=\frac{\beta}{y}, \beta\neq-\frac{1}{2},$ and $F(t,a,b)=\delta_1+\gamma_1a+\kappa_1 b.$ Then $$u_f(y)=\frac{\alpha}{1+2\beta}\left(\left(\frac{y}{\alpha}\right)^{1+2\beta}-1\right)\ \ \textmd{and} \ \ G_f^F(t,y,z)=\frac{\alpha^{2\beta}(\delta_1-\frac{\gamma_1\alpha}{1+2\beta})}{y^{2\beta}}+\frac{\gamma_1}{1+2\beta}y+\kappa_1 z.$$
}
\ \ \ \ \emph{(iv) Let $D=(0,\infty),$ $f(y)=-\frac{1}{2y},$ and $F(t,a,b)=\delta_1+\gamma_1a+\kappa_1 b.$ Then
$$u_f(y)=\alpha\ln\left(\frac{y}{\alpha}\right) \ \ \textmd{and} \ \ G_f^F(t,y,z)=\frac{\delta_1}{\alpha}y+\gamma_1 y\ln\left(\frac{y}{\alpha}\right)+\kappa_1 z.$$
}
\ \ \ \ \emph{(v) Let} $D=\textbf{R},$\emph{ $f(y)=\frac{\beta}{2},\ \beta\neq0,$ and $F(t,a,b)=\delta_1+\gamma_1a+\kappa_1 b.$ Then
$$u_f(y)=\frac{1}{\beta}(\exp(\beta(y-\alpha))-1)\ \ \textmd{and} \ \ G_f^F(t,y,z)=\frac{\delta_1-\frac{\gamma_1}{\beta}}{\exp(\beta(y-\alpha))}+\frac{\gamma_1}{\beta}+\kappa_1 z.$$
}

In this paper, we study the following RBSDE$(G_f^F+f(y)|z|^2,\xi,L_t)$:
$$\left\{
   \begin{array}{ll}
     Y_t=\xi+\int_t^{T}(G_f^F(s,Y_s,Z_s)+f(Y_s)|Z_s|^2)ds+K_T-K_t-\int_t^{T}Z_sdB_s,\ \ t\in[0,T],\\
     \forall t\in[0,T],\ \ Y_t\geq L_t,\\
     \int_0^T(Y_t-L_t)dK_t=0,
   \end{array}
 \right.
\eqno(2.2)$$
where $T$ is the terminal time, $\xi$ is the terminal variable, and $L_t$ is the lower obstacle.\\\\
\textbf{Definition 2.2} \  A solution of the RBSDE$(G_f^F(t,y,z)+f(y)|z|^2,\xi,L_t)$ is a triple of processes $(Y_t,Z_t,K_t)\in{\mathcal{C}}_D
\times{H}^2_d\times{\cal{A}},$ which satisfies $\int_0^{T}|G_f^F(s,Y_s,Z_s)+f(Y_s)|Z_s|^2)|ds<\infty$ and (2.2).\\

When the RBSDE$(G_f^F(t,y,z)+f(y)|z|^2,\xi,L_t)$ is not restricted by $L_t$, it becomes the standard BSDE$(g,\xi):$
$$Y_t=\xi +\int_t^T(G_f^F(s,Y_s,Z_s)+f(Y_s)|Z_s|^2)ds-\int_t^TZ_sdB_s,\ \ t\in[0,T].\eqno(2.3)$$
\textbf{Definition 2.3}  A solution of the BSDE$(G_f^F(t,y,z)+f(y)|z|^2,\xi)$ is a pair of processes $(Y_t,Z_t)\in{\mathcal{C}}_D
\times{H}^2_d,$ which satisfies $\int_0^{T}|G_f^F(s,Y_s,Z_s)+f(Y_s)|Z_s|^2|ds<\infty$ and (2.3).
%%%%%%%%%%%%%%%%%%%%%%%%%%%%%%%%%%%%%%%%%%%%%%%%%%%%%%%%%%%%%%%%%%%%%%%%%%%%%%%%%%%%%%%%%%%%%%%%%%%%%%%%%%%%%%%%%%%%%%%%%%%%%%%%%%%%%%%%%%%%%%%%%%%%%%%%%%%%%%%%%%
\section{Existence and uniqueness}
In this section, we study the existence and uniqueness of the solution to the RBSDE$(G_f^F+f(y)|z|^2,\xi,L_t)$. First, we provide a necessary condition for the existence.  \\\\
\textbf{Proposition 3.1}\ \emph{Let $(Y_t,Z_t,K_t)$ be a solution to the RBSDE$(f(y)|z|^2,\xi,L_t)$. If there exists a constant $\beta$ such that, for each $x\in D$, $u_f(x)\geq\beta$, then we have $u_f(\xi)\in L^1({\cal{F}}_T).$}
\\\\
\emph{Proof.} Applying Lemma A.2(ii) to $u_{f}(Y_t),$ and then by Lemma A.1(iii), we have
$$u_f(Y_t)=u_f(\xi)+\int_t^Tu'_f(Y_s)dK_s-\int_t^Tu'_f(Y_s)Z_sdB_s,\ \ t\in[0,T].$$
Set
$$\tau_n=\inf\left\{t\geq0,\int_0^t|u'_f(Y_s)|^2|Z_s|^2ds\geq n\right\}\wedge T,\ \ n\geq1.$$
By the two equalities above and the fact that $\int_0^tu'_f(Y_s)dK_s\geq0$, we have $u_f(Y_0)\geq E[u_f(Y_{\tau_n})].$ Then by the continuity of $Y_t$ and $u_f$, and Fatou's lemma, we have
$$\beta \leq E[u_f(\xi)]=E[\liminf_{n\rightarrow\infty}u_f(Y_{\tau_n})]\leq \liminf_{n\rightarrow\infty}E[u_f(Y_{\tau_n})]\leq u_f(Y_0),$$
which implies $u_f(\xi)\in L^1({\cal{F}}_T).$\ \ $\Box$\\

The following existence and uniqueness result plays an important role in this paper. \\\\
\textbf{Proposition 3.2} \emph{Let $e^{\gamma T}(\xi^++\delta T)\in L_{D}({\cal{F}}_T)$ with $\xi\in L^p({\cal{F}}_T)$, and $e^{\gamma t}(L^+_t+\delta t)\in{\mathcal{C}}_D$ with $L_t\in{\mathcal{S}}^p$. Then the RBSDE$(F,\xi,L_t)$ has a unique solution $(y_t,z_t,k_t)$ such that $y_t\in{\mathcal{S}}_D^p$. Moreover, we have $(z_t,k_t)\in{\cal{H}}^p\times{\mathcal{A}}^p$.}\\\\
\emph{Proof.} According to Bouchard et al. \cite[Theorem 3.1]{Bou}, the RBSDE$(F,\xi,L_t)$ has a solution $(y_t,z_t,k_t)$ such that $y_t\in{\mathcal{S}}^p$. If the RBSDE$(F,\xi,L_t)$ has another solution $(\tilde{y}_t,\tilde{z}_t,\tilde{k}_t)$ such that $\tilde{y}_t\in{\mathcal{S}}^p$, then by \cite[Proposition 2.1]{Bou} and a localization argument, we have $(y_t,z_t,k_t),(\tilde{y}_t,\tilde{z}_t,\tilde{k}_t)\in{\mathcal{S}}^p\times{\cal{H}}^p\times{\mathcal{A}}^p$. By \cite[Theorem 3.1]{Bou} again, we obtain $(y_t,z_t,k_t)=(\tilde{y}_t,\tilde{z}_t,\tilde{k}_t),$ $dt\times dP$-$a.e..$ In other words, the RBSDE$(F,\xi,L_t)$ has a unique solution $(y_t,z_t,k_t)$ such that $y_t\in{\mathcal{S}}^p$. Moreover, we have $(z_t,k_t)\in{\cal{H}}^p\times{\mathcal{A}}^p$.

We now consider the range of $y_t$. Using a classic linearization method, we have
$$y_t=\xi+\int_t^T(F(s,0,0)+a_sy_s+b_sz_s)ds+k_T-k_t-\int_t^Tz_sd{B}_s,\ \ t\in[0,T],\eqno(3.1)$$
where $$a_s=\frac{F\left(s,y_s,z_s\right)-F\left(s,0,z_s\right)}{y_s}1_{\{|y_s|>0\}}\ \ \ \textmd{and}\ \ \
b_s=\frac{(F\left(s,0,z_s\right)-F\left(s,0,0\right))z_s}{|z_s|^2}1_{\{|z_s|>0\}}.$$
Clearly, $|a_s|\leq \gamma, |b_s|\leq \kappa$. Let $Q$ be a probability measure such that $$\frac{dQ}{dP}=\exp\left\{\int_0^Tb_sdB_s-\frac{1}{2}\int_0^T|b_s|^2ds\right\}.$$
Then $X_t:=\exp\{\int_0^tb_sdB_s-\frac{1}{2}\int_0^t|b_s|^2ds\}$ is a solution to the the following linear SDE:
$$X_t=1+\int_0^tb_sX_sdB_s,\ \ t\in[0,T],$$
such that for each $r>1$, $X_t\in{\mathcal{S}}^r$ (see \cite[Theorem 4.4, page 61 and Lemma 2.3, page 93]{M}).
By Girsanov's theorem, $\bar{B}_t=B_t-\int_0^tb_sds$ is a standard Brownian motion under $Q.$ Then, (3.1) can be rewritten as
$$y_t=\xi+\int_t^T(F(s,0,0)+a_sy_s)ds+k_T-k_t-\int_t^Tz_sd\bar{B}_s,\ \ t\in[0,T].$$
By It\^{o}'s formula, we have
$$e^{\int_0^ta_rdr}y_t=e^{\int_0^Ta_rdr}\xi+\int_t^TF(s,0,0) e^{\int_0^sa_rdr}ds+\int_t^Te^{\int_0^sa_rdr}dk_s-\int_t^Te^{\int_0^sa_rdr}z_sd\bar{B}_s,\ \ t\in[0,T].\eqno(3.2)$$
Since $|a_s|\leq \gamma$ and for each $r>1$, $X_t\in{\mathcal{S}}^r$, by H\"{o}lder's inequality, we have, for $1<q<p,$
$$E_Q\left[\sup_{{t\in[0,T]}}|e^{\int_0^ta_rdr}y_t|^q\right]\leq e^{\gamma Tq}E\left[X_T\sup_{{t\in[0,T]}}|y_t|^q\right]\leq e^{\gamma Tq}E\left[|X_T|^\frac{p}{p-q}\right]^{^\frac{p-q }{p}}E\left[\sup_{{t\in[0,T]}}|y_t|^p\right]^{\frac{q}{p}}<\infty.$$
Similarly, we can also get $e^{\int_0^ta_rdr}z_t\in{\cal{H}}^q, \int_0^te^{\int_0^sa_rdr}dk_s\in{\mathcal{A}}^q$ and $e^{\int_0^ta_rdr}L_t\in{\mathcal{S}}^q$ under probability measure $Q.$ Thus, $(e^{\int_0^ta_rdr}y_t,e^{\int_0^ta_rdr}z_t,\int_0^te^{\int_0^sa_rdr}dk_s)\in{\mathcal{S}}^q
\times{\cal{H}}^q\times{\mathcal{A}}^q$ under probability measure $Q$ is a solution to the RBSDE$(F(t,0,0)e^{\int_0^ta_rdr},e^{\int_0^Ta_rdr}\xi,e^{\int_0^ta_rdr}L_t).$
Then, by the proof of \cite[Proposition 2.3]{E} and (3.2), we have
$$e^{\int_0^ta_rdr}y_t=\textrm{ess}\sup_{\tau\in{\cal{T}}_{t,T}}E_Q\left[\int_t^\tau F(s,0,0) e^{\int_0^sa_rdr}ds+e^{\int_0^\tau a_rdr}\eta_\tau|{\cal{F}}_t\right], \ \ \forall t\in[0,T],$$
where $\eta_\tau:=L_\tau1_{\{\tau<T\}}+\xi1_{\{\tau=T\}}$. Thus, we have
$$y_t=\textrm{ess}\sup_{\tau\in{\cal{T}}_{t,T}}E_Q\left[\int_t^\tau F(s,0,0) e^{\int_t^sa_rdr}ds+e^{\int_t^\tau a_rdr}\eta_\tau|{\cal{F}}_t\right], \ \ \forall t\in[0,T].\eqno(3.3)$$
Since $|F(s,0,0)|\leq \delta$ and $|a_s|\leq \gamma,$ we have
\begin{eqnarray*}
\ \ \ \   \ \ \ \  \ \ \ \ \ \   \ \ \ \  \ \ \  L_t\leq y_t&\leq& E_Q[\sup_{\tau\in{\cal{T}}_{t,T}}e^{\gamma(\tau-t)}(\delta (\tau-t)+\eta^+_\tau)|{\cal{F}}_t]\\
&\leq& E_Q[\sup_{\tau\in{\cal{T}}_{0,T}}e^{\gamma\tau}(\delta \tau+L^+_\tau)\vee e^{\gamma T}(\delta T+\xi^+)|{\cal{F}}_t], \ \ \forall t\in[0,T]. \ \ \ \ \ \ \ \ \ \ \ (3.4)
\end{eqnarray*}
We can prove that, for each $\eta\in L^1({\cal{F}}_T)$ and $c\in\textbf{R}$ such that $\eta<c$, we have $E[\eta|{\cal{F}}_t]<c$ for each $t\in[0,T]$. Thus, (3.4) leads to $y_t\in{\mathcal{S}}_D$. The proof is complete. \ \ $\Box$\\

We have the following existence and uniqueness result for the RBSDE$(G_f^F+f(y)|z|^2,\xi,L_t)$.\\\\
\textbf{Theorem 3.3}  \emph{Let $e^{\gamma T}(u_f(\xi)\vee0+\delta T)\in L_V({\cal{F}}_T)$ with $u_f(\xi)\in L^p({\cal{F}}_T)$, and $e^{\gamma t}(u_f(L_t)\vee0+\delta t)\in{\mathcal{C}}_{V}$ with $u_f(L_t)\in{\mathcal{S}}^p$. Then the RBSDE$(G_f^F+f(y)|z|^2,\xi,L_t)$ has a unique solution $(Y_t,Z_t,K_t)$ such that $u_f({Y}_t)\in{\mathcal{S}}^p$. Moreover, we have}

\emph{(i) $(Y_t,Z_t,K_t)\in{\mathcal{S}}^p\times{\cal{H}}^p\times{\mathcal{A}}^p,$ when $f$ is integrable on $D$;}

\emph{(ii) $(Y_t,Z_t,K_t)\in{\mathcal{S}}^p\times{\cal{H}}^{2p}\times{\mathcal{A}}^p,$ when there exist constants $c\in D$ and $\beta>0$ such that $L_t\geq c, dt\times dP-a.e.$ and $f\geq\beta, a.e..$}\\\\
\emph{Proof.} By Proposition 3.2, the RBSDE$(F,u_f(\xi),u_f(L_t))$ has a unique solution $(y_t,z_t,k_t)$ such that $y_t\in{\mathcal{S}}_V^p.$ Moreover, we have $(z_t,k_t)\in{\cal{H}}^p\times{\mathcal{A}}^p$. In view of Lemma A.1(iv), we can apply Lemma A.2(ii) to $u_f^{-1}(y_t),$ and then, by setting
$$(Y_t,Z_t,K_t):=\left(u_f^{-1}(y_t),\frac{z_t}{u_f'(u_f^{-1}(y_t))}, \int_0^t\frac{1}{u_f'(u_f^{-1}(y_s))}dk_s\right),\ \ t\in[0,T],\eqno(3.5)$$
and using Lemma A.1(iii), we get that $(Y_t,Z_t,K_t)$ is a solution to the RBSDE$(G_f^F+f(y)|z|^2,\xi,L)$ such that $u_f({Y}_t)\in{\mathcal{S}}^p.$ We now prove the uniqueness. Since $u_f(x)\in W^2_{1,loc}(D)$ (see Lemma A.1(ii)), for a solution $(\bar{Y}_t,\bar{Z}_t,\bar{K}_t)$ of RBSDE$(G_f^F+f(y)|z|^2,\xi,L)$ such that $u_f(\bar{Y}_t)\in{\mathcal{S}}^p,$
applying Lemma A.2(ii) to $u_f(\bar{Y}_t)$, and then by Lemma A.1(iii), we have
$$u_f(\bar{Y}_t)=u_f(\xi)+\int_t^TF(s,u_f(\bar{Y}_t),u'_f({\bar{Y}}_s)\bar{Z}_s)ds+\int_t^Tu_f'(\bar{Y}_s)d\bar{K}_s-\int_t^Tu'_f({\bar{Y}}_s)\bar{Z}_sdB_s,\ \ t\in[0,T],$$ which means that $(u_f(\bar{Y}_t),u'_f(\bar{Y}_t)\bar{Z}_t,\int_0^tu_f'(\bar{Y}_s)d\bar{K}_s)$ is a solution to the RBSDE$(F,u_f(\xi),u_f(L_t))$ such that $u_f(\bar{Y}_t)\in{\mathcal{S}}^p.$ Moreover, from (3.5), it follows that $(u_f(Y_t),u'_f(Y_t)\bar{Z}_t,\int_0^tu_f'(Y_s)d\bar{K}_s)$ is a unique solution to the RBSDE$(F,u_f(\xi),u_f(L_t))$ such that $u_f({Y}_t)\in{\mathcal{S}}^p$. Then, by Lemma A.1(ii), we obtain the uniqueness.

Proof of (i): Since $f$ is integrable on $D$, by \cite[Lemma A.1(j)]{B17}, there exist two positive constants $c_1$ and $c_2$ such that, for each $x,y\in D$, we have
$c_1|x-y|\leq |u_f(x)-u_f(y)|\leq c_2|x-y|.$
Then, by (3.5), Lemma A.1(ii), and the fact that $(y_t,z_t,k_t)\in{\mathcal{S}}^p
\times{\cal{H}}^p\times{\mathcal{A}}^p$, we get (i).

Proof of (ii): Since $f>0, a.e.$, by Lemma A.1(v), we have
$$u_f(Y_t)\geq u_0(Y_t)=Y_t-\alpha\geq L_t-\alpha,$$
which together with $L_t\geq c, dt\times dP-a.e.,$ implies $Y_t\in{\mathcal{S}}^p$ and $Y_t\geq c, dt\times dP-a.e.$

For $n\geq1,$ we define the following stopping time
$$\sigma_n=\inf\left\{t\geq0,\int_0^tf(Y_s)|Z_s|^2ds\geq n\right\}\wedge T.$$
Since $Y_t\geq c,\ dt\times dP-a.e.$, by the assumption $f\geq\beta, a.e.,$ (2.1), (2.2), and the fact that $K_s\in\cal{A}$, for any stopping time $\tau\leq\sigma_1$, we have
\begin{eqnarray*} \ \ \ \ \ \ \ \ \ \ \beta\int_\tau^{\sigma_n}|Z_s|^2ds&\leq&\int_\tau^{\sigma_n}f(Y_s)|Z_s|^2ds
\\&=&Y_\tau-Y_{\sigma_n}-\int_\tau^{\sigma_n}G_f^F(s,Y_s,Z_s)ds-(K_{\sigma_n}-K_\tau)+\int_\tau^{\sigma_n}Z_sdB_s
\\&\leq&|Y_\tau|+|Y_{\sigma_n}|+\int_\tau^{\sigma_n}C(1+|Y_s|+|Z_s|)ds+\int_\tau^{\sigma_n}Z_sdB_s.\ \ \ \ \ \ \ \ \ \ \ \ \ \ \ \ (3.6)
\end{eqnarray*}
where $C$ is a constant depending only on $c, u_f'(c), u_f(c), \delta, \gamma,$ and $\kappa.$ By Jensen's inequality, (3.6), and BDG inequality, we have
\begin{eqnarray*}
\ \ \ \ \beta^p\left(E\left[\left(\int_{\tau}^{\sigma_n}|Z_s|^2ds\right)^{\frac{p}{2}}\right]\right)^2&\leq& E\left[\left(\beta\int_{\tau}^{\sigma_n}|Z_s|^2ds\right)^p\right]\\&\leq&E\left[\left(\int_\tau^{\sigma_n}f(Y_s)|Z_s|^2ds\right)^p\right]\\&\leq&C_1\left(1
+E\left[\left(\int_{\tau}^{\sigma_n}|Z_s|^2ds\right)^{\frac{p}{2}}\right]+E\left(\left|\int_{\tau}^{\sigma_n}Z_sdB_s\right|^p\right)\right)\\
&\leq&C_2\left(1+E\left[\left(\int_{\tau}^{\sigma_n}|Z_s|^2ds\right)^{\frac{p}{2}}\right]\right),\ \ \ \ \ \ \ \ \ \ \ \ \ \ \ \ \ \ \ \ \ \ \ \ \ \ \ \ \ (3.7)
\end{eqnarray*}
where $C_1$ and $C_2$ are two positive constants depending only on $E[\sup_{t\in[0,T]}|Y_t|^p]$, $T, C,$ and $p$. Then by solving the quadratic inequality (3.7) with $E[(\int_{\tau}^{\sigma_n}|Z_s|^2ds)^{\frac{p}{2}}]$ as the unknown variable, we get $E[(\int_{\tau}^{\sigma_n}|Z_s|^2ds)^{\frac{p}{2}}]<C_3,$ where $C_3>0$ is a constant depending only on $C_2$ and $\beta$. By plugging this inequality into (3.7), we get $E\left[\left(\int_{\tau}^{\sigma_n}|Z_s|^2ds\right)^p\right]<C_4,$ where $C_4>0$ is a constant depending only on $C_2,C_3,$ and $\beta$. Fatou's Lemma then gives
$$E\left[\left(\int_0^{T}|Z_s|^2ds\right)^p\right]<\infty. \eqno(3.8)$$
Since $f>0, a.e.,$ by (2.1) and (2.2), we have
\begin{eqnarray*}
0\leq K_{T}&\leq&Y_0-\xi-\int_0^{T}G_f^F(s,Y_s,Z_s)ds+\int_0^{T}Z_sdB_s\\
&\leq&|Y_0|+|\xi|+\int_0^{T}C(1+|Y_s|+|Z_s|)ds+\left|\int_0^{T}Z_sdB_s\right|.
\end{eqnarray*}
where $C$ is the constant in (3.6). Then by the fact that $Y_t\in{\mathcal{S}}^p$, BDG inequality, and (3.8), we get $K_t\in{\mathcal{A}}^p$. Thus, (ii) holds. The proof is complete.\ \ $\Box$
\\

From Theorem 3.3, we get the following Corollary 3.4 directly.\\\\
\textbf{Corollary 3.4}  \emph{Let $V=(b,\infty)$ with $b\geq-\infty$. Let $u_f(\xi)\in L^p({\cal{F}}_T)$ and $u_f(L_t)\in{\mathcal{S}}^p$. Then the RBSDE$(G_f^F+f(y)|z|^2,\xi,L_t)$ has a unique solution $(Y_t,Z_t,K_t)$ such that $u_f({Y}_t)\in{\mathcal{S}}^p$.}\\\\
\textbf{Corollary 3.5}  \emph{Let $\delta=\gamma=0$ in Assumption (A). Let $u_f(\xi)\in L^p({\cal{F}}_T)$ and $u_f(L_t)\in{\mathcal{S}}^p$. Then the RBSDE$(G_f^F+f(y)|z|^2,\xi,L_t)$ has a unique solution $(Y_t,Z_t,K_t)$ such that $u_f({Y}_t)\in{\mathcal{S}}^p$. In particular, for each $t\in[0,T]$, we have}
$$u_f(Y_t)=\textrm{ess}\sup_{\tau\in{\cal{T}}_{t,T}} E_Q[u_f(\eta_\tau)|{\cal{F}}_t],\eqno(3.9)$$
\emph{where $\eta_\tau:=L_\tau1_{\{\tau<T\}}+\xi1_{\{\tau=T\}}$ and $Q$ is a probability measure equivalent to $P$.}\\\\
\emph{Proof.} Since $\delta=\gamma=0$ and $0\in V$, by Theorem 3.3, we get that the RBSDE$(G_f^F+f(y)|z|^2,\xi,L_t)$ has a unique solution $(Y_t,Z_t,K_t)$ such that $u_f({Y}_t)\in{\mathcal{S}}^p$. By (3.3) and (3.5), we deduce that, for each $t\in[0,T]$,
$$u_f(Y_t)=\textrm{ess}\sup_{\tau\in{\cal{T}}_{t,T}}E_Q\left[\int_t^\tau F(s,0,0) e^{\int_t^sa_rdr}ds+e^{\int_t^\tau a_rdr}u_f(\eta_\tau)|{\cal{F}}_t\right],$$
where $\eta_\tau:=L_\tau1_{\{\tau<T\}}+\xi1_{\{\tau=T\}}$, $Q$ is a probability measure equivalent to $P$, and $a_t$ is a process such that $|a_t|\leq\gamma$. Since $\delta=\gamma=0$, we have $F(s,0,0)=0$ and $a_s=0$. The proof is complete. \ \ $\Box$ \\\\
\textbf{Corollary 3.6}  \emph{Let $\delta=\gamma=0$ in Assumption (A) and $u_f(\xi)\in L^p({\cal{F}}_T)$. Then the BSDE$(G_f^F+f(y)|z|^2,\xi)$ has a unique solution $(Y_t,Z_t)$ such that $u_f({Y}_t)\in{\mathcal{S}}^p$. In particular, for each $t\in[0,T]$, we have}
$$u_f(Y_t)=E_Q[u_f(\xi)|{\cal{F}}_t],\eqno(3.10)$$
\emph{where $Q$ is a probability measure equivalent to $P$.}\\\\
\emph{Proof.} By \cite[Theorem 4.2 and Lemma 3.1]{BD}, we deduce
that the BSDE$(F, u_f(\xi))$ has a unique solution $(y_t, z_t)$ such that $y_t\in{\cal{S}}^p.$ Since $\delta=\gamma=0$, from (3.2), we get that
$$y_t=E_Q\left[u_f(\xi)|{\cal{F}}_t\right],$$
where $Q$ is a probability measure equivalent to $P$. This implies that $y_t\in{\cal{S}}_V^p.$ Then applying Lemma A.2(ii) to $u_f^{-1}(y_t),$ and by a similar argument as in the proof of Theorem 3.3, we deduce that the BSDE$(G_f^F+f(y)|z|^2,\xi)$ has a unique solution $(Y_t,Z_t)$ such that $u_f({Y}_t)=y_t$. The proof is complete. \ \ $\Box$ \\\\
\textbf{Remark 3.7} \emph{According to the definition of $F$, Corollary 3.6 is an extension of \cite[Proposition 3.3]{Zheng1}, which studied the BSDE$(K|z|+f(y)|z|^2,\xi)$, where $K$ is a constant. In Corollary 3.5 and 3.6, if we further assume that $\kappa=0$ in Assumption (A), then the probability measure $Q$ in (3.9) and (3.10) are both the probability measure $P$.}\\\\
\textbf{Example 3.8}\ \ We show some typical cases of Corollary 3.4.

(i) For $f(y)=\frac{\beta}{y}$ with $\beta>-\frac{1}{2}$, if $\xi^{1+2\beta}\in L^p({\cal{F}}_T)$ and $L_t^{1+2\beta}\in{\mathcal{S}}^p,$ then the RBSDE$(G_f^F+f(y)|z|^2,\xi,L_t)$ has a unique solution $(Y_t,Z_t,K_t)$ such that $Y_t^{1+2\beta}\in{\mathcal{S}}^p.$

(ii) For $f(y)=-\frac{1}{2y}$, if $\ln(\xi)\in L^p({\cal{F}}_T)$ and $\ln(L_t)\in{\mathcal{S}}^p,$ then the RBSDE$(G_f^F+f(y)|z|^2,\xi,L_t)$ has a unique solution $(Y_t,Z_t,K_t)$ such that $\ln(Y_t)\in{\mathcal{S}}^p.$

(iii) For $f(y)=\frac{\beta}{2}$ with $\beta>0$, if $\exp(\beta\xi)\in L^p({\cal{F}}_T)$ and $\exp(\beta L_t)\in{\mathcal{S}}^p,$ then the RBSDE$(G_f^F+f(y)|z|^2,\xi,L_t)$ has a unique solution $(Y_t,Z_t,K_t)$ such that $\exp(\beta Y_t)\in{\mathcal{S}}^p.$

(iv) For $f(y)=\frac{\beta}{2}1_{\{y\geq0\}}$ with $\beta>0$, if $\xi^-, \exp(\beta\xi^+)\in L^p({\cal{F}}_T)$ and $L_t^-, \exp(\beta L_t^+)\in{\mathcal{S}}^p,$ then the RBSDE$(G_f^F+f(y)|z|^2,\xi,L_t)$ has a unique solution $(Y_t,Z_t,K_t)$ such that $Y_t^-, u_f(\beta Y_t^+)\in{\mathcal{S}}^p.$\\

The following Example 3.9 shows that the quadratic BSDE with a bounded terminal variable may have no bounded solution, even when $D=\textbf{R}$, $f(y)$ is a constant, and the corresponding RBSDE has a bounded solution. \\\\
\textbf{Example 3.9}\ \ We consider the Example 2.1(v) with $\beta=1, \delta_1=0, e^{\gamma_1 T}>2$ and $\alpha=0$. Corollary 3.4 shows that when $u_f(L_t)\in{\mathcal{S}}^p$ and $L_T\leq\ln(\frac{1}{2})$, the RBSDE$(G_f^F+\frac{|z|^2}{2},\ln(\frac{1}{2}), L_t)$ has a unique solution $(Y_t,Z_t,K_t)$ such that $u_f( Y_t)\in{\mathcal{S}}^p.$ Moreover, if $L_t$ is bounded, by (3.4) and (3.5), we deduce that $Y_t$ is bounded. However, we will see that the BSDE$(G_f^F+\frac{|z|^2}{2},\ln(\frac{1}{2}))$ has no solution $(Y_t,Z_t)$ such that $u_f(Y_t)\in{\mathcal{S}}^p$. In fact, by \cite[Theorem 4.2 and Lemma 3.1]{BD}, we deduce that the BSDE$(F,-\frac{1}{2})$ has a unique solution $(y_t,z_t)$ such that $y_t\in{\mathcal{S}}^p$. Since $u_f(x)=\exp(x)-1$ and $V=(-1,\infty)$, from (3.2), we get that
    $$y_0=E_Q\left[-\frac{1}{2}e^{\int_0^T \gamma_1dr}\right]=-\frac{1}{2}e^{\gamma_1 T}\notin V,\eqno(3.11)$$
where $Q$ is a probability measure equivalent to $P$. If we assume that the BSDE$(G_f^F+\frac{|z|^2}{2},\ln(\frac{1}{2}))$ has a solution $(Y_t,Z_t)$ such that $u_f(Y_t)\in{\mathcal{S}}^p$, then by applying It\^{o}'s formula to $u_f(Y_t)$, we get that $(u_f(Y_t),u'_f(Y_t)Z_t)$ is a solution to the BSDE$(F,-\frac{1}{2})$. This together with the uniqueness of the solution to the BSDE$(F,-\frac{1}{2})$ gives $y_t=u_f(Y_t)\in V$, which contradicts (3.11).\\

We now establish some comparison results for the RBSDE$({G}_f^F+f(y)|z|^2,\xi)$.\\\\
\textbf{Proposition 3.10} \emph{Let $F_1$ and $F_2$ satisfy Assumption (A) such that $F_1(\cdot)\geq F_2(\cdot)$. Let $\xi_1, \xi_1\in L_D({\cal{F}}_T)$ such that $\xi_1\geq\xi_2\geq L_T$ and $u_f(\xi_1), u_f(\xi_2)\in L^p({\cal{F}}_T)$. Let $(Y^i_t,Z^i_t,K^i_t)$ be the solution of the RBSDE$({G}_f^{F_i}+f(y)|z|^2,\xi_i,L)$ satisfying $u_{f}(Y^i_t)\in{\mathcal{S}}^p, i=1,2.$ Then for each $t\in[0,T]$, $Y^1_t\geq Y^2_t$.}\\\\
\emph{Proof.} By the proof of Theorem 3.3, we can get that the RBSDE$(F_i,u_f(\xi_i),u_f(L_t))$ has a unique solution $(y^i_t,z^i_t,k^i_t)$ such that $y_t^i=u_{f}(Y^i_t),$ $i=1,2.$ By Lemma A.1(ii) and \cite[Theorem 4.1]{E}, we get that for each $t\in[0,T]$, $u_{f}(Y_t^1)\geq u_{f}(Y^2_t),$ which gives $Y^1_t\geq Y^2_t$.\ \ $\Box$\\

We assume two semimartingales:
$$Y^1_t=Y^1_T+\int_t^Th_1(s)ds+A^1_T-A^1_t-\int_t^TZ^1_sdB_s,\ \ t\in[0,T],$$
and
$$Y^2_t=Y^2_T+\int_t^Th_2(s)ds-A^2_T+A^2_t-\int_t^TZ^2_sdB_s,\ \ t\in[0,T],$$
where $Y_t^i\in{\cal{C}}_D$, $A^i_t\in {\cal{A}},$ $Z^i_t\in H^2$ and $h_i(t)$ is a progressively measurable process satisfying $\int_0^T|h_i(t)|dt<\infty, i=1,2.$\\\\
\textbf{Proposition 3.11} \emph{Let $u'_{f}(Y^i_t)Z^i_t\in{\cal{H}}^p$ and $u'_f(Y_t^i)h_i(t)-\frac{1}{2}u''_f(Y_t^i)|Z_t^i|^2\in{\cal{H}}^p$, $i=1,2$. Let $(Y_t,Z_t)$ be the solution of the BSDE$({G}_f^F+f(y)|z|^2,\xi)$ satisfying $u_{f}(Y_t)\in{\mathcal{S}}^p.$}

\emph{(i) Assume that $\xi\leq Y_T^1$ and $G_f^F(t,Y_t^1,Z_t^1)+f(Y_t^1)|Z_t^1|^2\leq h_1(t),\ dt\times dP$-a.e.. Then for each $t\in[0,T]$, $Y_t\leq Y_t^1.$ Moreover, if $Y_t=Y_t^1$, then $\xi=Y^1_T,$ $A_T^1=A_t^1,$ and $G_f^F(s,Y_s^1,Z_s^1)+f(Y_s^1)|Z_s^1|^2=h_1(s),\ dt\times dP$-a.e., on $[t,T]\times\Omega$;}

\emph{(ii) Assume that $\xi\geq Y_T^2$ and $G_f^F(t,Y_t^2,Z_t^2)+f(Y_t^2)|Z_t^2|^2\geq h_2(t),\ dt\times dP$-$a.e..$ Then for each $t\in[0,T]$, $Y_t\geq Y_t^2.$ Moreover, if $Y_t=Y_t^2,$ then $\xi=Y_T^2,$ $A_T^2=A_t^2,$ and $G_f^F(s,Y_s^2,Z_s^2)+f(Y_s^2)|Z_s^2|^2=h_2(s),\ dt\times dP$-$a.e.,$ on $[t,T]\times\Omega$.}\\\\
\emph{Proof.} Proof of (i): Applying Lemma A.2(ii) to $u_f(Y_t)$, we get that $(u_{f}(Y_t),u'_{f}(Y_t)Z_t)$ is a solution to the BSDE$(F,u_f(\xi))$. Since $u_{f}(Y_t)\in{\mathcal{S}}^p$, by \cite[Lemma 3.1]{BD}, we have $u'_{f}(Y_t)Z_t\in{\cal{H}}^p$. Applying Lemma A.2(ii) to $u_f(Y^1_t)$, we have
\begin{align*}
  u_{f}(Y^1_t)=u_f(Y^1_T)+\int_t^T(u'_f(Y_s^1)h_1(s)&-\frac{1}{2}u''_f(Y_s^1)|Z_s^1|^2)ds\\
&+\int_t^Tu'_f(Y_s^1)dA_s^1-\int_t^Tu'_{f}(Y^1_t)Z^1_sdB_s,\ \ \ t\in[0,T].
\end{align*}
Since for each $t\in[0,T]$,
\begin{align*}
F(t,u_{f}(Y^1_t),u'_{f}(Y^1_t)Z^1_t)&=G_f^F(t,Y^1_t,Z^1_t)u'_f(Y^1_t),\ \ \ (\textrm{by the definition of}\ G_f^F)\\
&\leq({h}_1(t)-f(Y_t^1)|Z_t^1|^2)u'_f(Y^1_t)\\
&=u'_f(Y_t^1)h_1(t)-\frac{1}{2}u''_f(Y_t^1)|Z_t^1|^2,\ \ \ (\textrm{by Lemma A.1(iii)}),
\end{align*}
it follows from the proof of \cite[Theorem 2.2]{EPQ}, and Lemma A.1(ii) that (i) holds.

Proof of (ii):  The proof is similar to (i), so it is not given explicitly. \ $\Box$\\\\
\textbf{Remark 3.12} \emph{If $f=0,$ then Proposition 3.10 becomes a comparison theorem for RBSDEs with Lipschitz continuous generators, which was studied by \cite[Theorem 4.1]{E}. If $G_f^F=0$, then Proposition 3.11 becomes the comparison theorem for the BSDE$(f(y)|z|^2,\xi)$, which was studied by \cite[Proposition 3.2]{B17}, \cite[Proposition 2.3]{B19}, and \cite[Proposition 4.1, 4.3]{Zheng1}.}
%%%%%%%%%%%%%%%%%%%%%%%%%%%%%%%%%%%%%%%%%%%%%%%%%%%%%%%%%%%%%%%%%%%%%%%%%%%%%%%%%%%%%%%%%%%%%%%%%%%%%%%%%%%%%%%%%%%%%%%%%%%%%%%%%%%%%%%%%%%%%%%%%%%%%%%%%%%%%%%%%%
\section{Applications}
\subsection{An obstacle problem for PDEs with singular coefficients}
In this subsection, we study the following obstacle problem for a quadratic PDE:
$$\left\{
    \begin{array}{ll}
      \min\{v(t,x)-h(t,x),-\partial_tv(t,x)-{\cal{L}}v(t,x)-\widetilde{G}(t,v(t,x),\sigma^*\nabla_xv(t,x))\}=0,\ (t,x)\in[0,T)\times \textbf{R}^d;\\
      v(T,x)=\psi(x),\ \ x\in \textbf{R}^d,
    \end{array}
  \right.\eqno(4.1)
$$
where $\widetilde{G}(t,y,z)=G_f^F(t,y,z)+f(y)|z|^2$ with $\delta=\gamma=0$, $\psi(x):\textbf{R}^d\mapsto D,$ $h(t,x):[0,T]\times \textbf{R}^d\mapsto D,$ and ${\cal{L}}$ is the infinitesimal generator of the solution $X_s^{t,x}$ of the SDE:	
$$X_s^{t,x}=x+\int_t^sb(r,X_r^{t,x})dr+\int_t^s\sigma(r,X_r^{t,x})dB_r,\ x\in \textbf{R}^d,\ s\in[t,T],$$
where $b:[0,T]\times \textbf{R}^d\mapsto \textbf{R}^d,$ $\sigma:[0,T]\times \textbf{R}^d\mapsto \textbf{R}^{d\times d}.$ The operator ${\cal{L}}$ is given by
$${\cal{L}}:=\frac{1}{2}\sum_{i,j=1}^d(\sigma\sigma^*)_{i,j}(s,x)\frac{\partial^2}{\partial_{x_i}\partial_{x_j}}
+\sum_{i=1}^db_{i}(s,x)\frac{\partial}{\partial_{x_i}}.$$
\textbf{Definition 4.1}\ \ A function $v(t,x)\in C_D([0,T]\times \textbf{R}^d)$\footnote{$v(t,x)\in C_D([0,T]\times \textbf{R}^d)$ means that $v(t,x)\in C([0,T]\times \textbf{R}^d)$ takes values in $D$.} is called a viscosity subsolution of (4.1), if $v(T,\cdot)\leq \psi(\cdot)$ and for each $(t,x,\phi)\in[0,T]\times \textbf{R}^d\times C^{1,2}_D([0,T]\times \textbf{R}^d)$ such that $\phi(t,x)=v(t,x)$ and $(t,x)$ is a local minimum point of $\phi-v$, we have
$$\min\{\phi(t,x)-h(t,x),-\partial_t\phi(t,x)-{\cal{L}}\phi(t,x)-\widetilde{G}(t,\phi(t,x),\sigma^*\nabla_x\phi(t,x))\}\leq0.$$

A function $v(t,x)\in C_D([0,T]\times \textbf{R}^d)$ is called a viscosity supersolution of (4.1), if $v(T,\cdot)\geq \psi(\cdot)$ and for each $(t,x,\phi)\in[0,T]\times \textbf{R}^d\times C^{1,2}_D([0,T]\times \textbf{R}^d)$ such that $\phi(t,x)=v(t,x)$ and $(t,x)$ is a local maximum point of $\phi-v$, we have
$$\min\{\phi(t,x)-h(t,x),-\partial_t\phi(t,x)-{\cal{L}}\phi(t,x)-\widetilde{G}(t,\phi(t,x),\sigma^*\nabla_x\phi(t,x))\}\geq0.$$

A function $v(t,x)\in C_D([0,T]\times \textbf{R}^d)$ is called a viscosity solution of (4.1), if it is a viscosity subsolution and a viscosity supersolution of (4.1).\\\\
\textbf{Assumption (B)} (i) $f$ is continuous.

(ii) $\psi(\cdot)$ is continuous such that $\psi(\cdot)\geq h(T,\cdot)$, and $u_f(\psi(\cdot))$ has polynomial growth.

(iii) $u_f(h(\cdot,\cdot))$ is continuous and $u_f(h(t,\cdot))$ has polynomial growth (uniformly in $t$).

(iv) $b(t,\cdot)$ and $\sigma(t,\cdot)$ are both Lipschitz continuous with linear growth (uniformly in $t$).\\

Let Assumption (B) hold. For $t\in[0,T)$, let $({\cal{F}}^t_s)_{t\leq s\leq T}$ be
the natural filtration generated by ${{(B_s-B_t)}_{s\geq t}}$, augmented
by the $\mathit{P}$-null sets of ${\cal{F}}$. Then by \cite[Theorem 4.4, page 61 and Lemma 2.3, page 93]{M}), for $(t,x)\in[0,T)\times \textbf{R}^d$, we have $u_f(\psi(X_T^{t,x}))\in L^2({\cal{F}}^t_T)$ and $u_f(h(s,X_s^{t,x}))\in{\mathcal{S}}^2$ ($({\cal{F}}^t_s)$-progressively measurable). It follows from Corollary 3.5 that the following Markovian RBSDE:
$$\left\{
   \begin{array}{ll}
     Y_s^{t,x}=\psi(X_T^{t,x})+\int_s^T\widetilde{G}(t,Y_r^{t,x},Z_r^{t,x})dr+K_T^{t,x}-K_s^{t,x}-\int_s^TZ_r^{t,x}dB_r, \ \ s\in[t,T],\\
     Y_s^{t,x}\geq h(s,X^{t,x}_s),\ \ s\in[t,T],\\
     \int_t^T(Y_s^{t,x}-h(s,X^{t,x}_s))dK_s^{t,x}=0,
   \end{array}
 \right.\eqno(4.2)$$
has a unique $({\cal{F}}^t_s)$-progressively measurable solution $(Y_s^{t,x},Z_s^{t,x},K_s^{t,x})$ such that $u_f(Y_s^{t,x})\in{\mathcal{S}}^2$. In addition, we have $u'_f(Y_s^{t,x})Z_s^{t,x}\in{\mathcal{H}}^2$. Set $v(t,x):=Y_t^{t,x}.$ We have the following proposition:\\\\
\textbf{Proposition 4.2}\ \ $v(t,x)$ is a viscosity solution of the obstacle problem (4.1). \\\\
\emph{Proof.} By (3.5), for each $(t,x)\in[0,T)\times \textbf{R}^d$, the RBSDE$(F,u_f(\psi(X^{t,x}_T)),u_f(h(s,X^{t,x}_s)))$ admits a unique $({\cal{F}}^t_s)$-progressively measurable solution $(y_s^{t,x},z_s^{t,x},k_s^{t,x})$ such that $y_s^{t,x}=u_f(Y_s^{t,x}).$ Thus, we have $u_f(v(t,x))=u_f(Y_t^{t,x})=y_t^{t,x}.$ It follows from \cite[Lemma 8.4]{E} and Lemma A.1(ii) that $v(t,x)$ is continuous in $(t,x)$. By (6.4.9) in \cite{Zhang}, we further get that, for $s\in[t,T],$
$$Y^{t,x}_s=v(s,X^{t,x}_s).\eqno(4.3)$$

Step 1. We will show that $v(t,x)$ is a viscosity supersolution of (4.1). Suppose that $v(t,x)$ is not a viscosity supersolution of (4.1). This means that there exists $(t,x,\phi)\in[0,T]\times \textbf{R}^d\times C^{1,2}_D([0,T]\times \textbf{R}^d)$ satisfying the condition that $\phi(t,x)=v(t,x)$ and that $(t,x)$ is a local maximum point of $\phi-v$, such that
$$\frac{\partial\phi}{\partial t}(t,x)+{\cal{L}}\phi(t,x)+\widetilde{G}(t,\phi(t,x),\sigma^*\nabla_x\phi(t,x))>0.$$
By continuity, there exist $\beta\in(0,T-t],$ $c>0,$ and $C>0$ such that for each $s\in[t,t+\beta]$ and $y\in[x-c,x+c]$, we have
$$v(s,y)\geq\phi(s,y)\ \ \textmd{and}\ \ \frac{\partial\phi}{\partial t}(s,y)+{\cal{L}}\phi(s,y)+\widetilde{G}(s,\phi(s,y),\sigma^*\nabla_x\phi(s,y))\geq C.\eqno(4.4)$$
We define a stopping time $\tau=\inf\{s\geq t; |X_s^{t,x}-x|\geq c\}\wedge(t+\beta).$ Thus, $t<\tau \leq t+\beta$, and the ranges of $v({s\wedge\tau},X_{s\wedge\tau}^{t,x})_{s\in[t,t+\beta]}$, $h({s\wedge\tau},X_{s\wedge\tau}^{t,x})_{s\in[t,t+\beta]}$ and $\phi({s\wedge\tau},X_{s\wedge\tau}^{t,x})_{s\in[t,t+\beta]}$ are all included in a closed subset of $D$. By (4.2) and (4.3), we get that the following RBSDE:
$$\bar{Y}_s=v(\tau,X_\tau^{t,x})+\int_s^{t+\beta} 1_{\{r\leq\tau\}}\widetilde{G}(r,\bar{Y}_r,\bar{Z}_r)dr+\bar{K}_{\tau}^{t,x}-\bar{K}_s^{t,x}-\int_s^{t+\beta}\bar{Z}_rdB_r,\ \ s\in[t,t+\beta],\eqno(4.5)$$
has a solution $(\bar{Y}_s,\bar{Z}_s,\bar{K}_s)=({Y}_{s\wedge\tau}^{t,x},1_{\{s\leq\tau\}}{Z}_s^{t,x},{K}_{s\wedge\tau}^{t,x})\in{\mathcal{S}}_D^\infty
\times{\cal{H}}^2\times{\mathcal{A}}$, where $\bar{Z}_s\in{\cal{H}}^2$ is due to $u'_f(Y_s^{t,x})Z_s^{t,x}\in{\mathcal{H}}^2$ and ${Y}_{s\wedge\tau}^{t,x}\in{\mathcal{S}}_D^\infty$. By Corollary 3.6, the BSDE$(\widetilde{G},v(\tau,X_\tau^{t,x}))$ with terminal time $\tau$ has a unique solution $(\hat{Y}_s,\hat{Z}_s)$ on $[t,\tau]$ such that $\hat{Y}_s\in{\mathcal{S}}_D^\infty$. Since $(\bar{Y}_s,\bar{Z}_s,\bar{K}_s)\in{\mathcal{S}}_D^\infty
\times{\cal{H}}^2\times{\mathcal{A}}$ and $F$ has a linear growth, we have $u'_{f}(\bar{Y}_t)\bar{Z}_t\in{\cal{H}}^2$ and
$$u'_f(\bar{Y}_t)\widetilde{G}(t,\bar{Y}_t,\bar{Z}_t)-\frac{1}{2}u''_f(\bar{Y}_t)|\bar{Z}_t|^2=F(t,u_{f}(\bar{Y}_t),u'_{f}(\bar{Y}_t)\bar{Z}_t)\in{\cal{H}}^2.$$
Then by (4.5) and Proposition 3.11(i), we get ${\bar{Y}}_t\geq\hat{Y}_t$.

Applying It\^{o}'s formula to $\phi(s,X_s^{t,x})$ for $s\in[t,\tau]$, we get that the following BSDE:
$${\tilde{Y}}_s=\phi(\tau,X_\tau^{t,x})+\int_s^{t+\beta} -1_{\{r\leq\tau\}}\left(\frac{\partial\phi}{\partial t}(r,X_r^{t,x})+{\cal{L}}\phi(r,X_r^{t,x})\right)dr-\int_s^{t+\beta}{\tilde{Z}}_rdB_r,\ \ s\in[t,t+\beta],\eqno(4.6)$$
has a solution $({\tilde{Y}}_s,{\tilde{Z}}_s)=(\phi({s\wedge\tau},X_{s\wedge\tau}^{t,x}), 1_{\{s\leq\tau\}}\sigma^*\nabla_x\phi(s,X_s^{t,x}))\in{\mathcal{S}}_D^\infty
\times {\cal{H}}^\infty.$ By (4.4), we have $v(\tau,X_\tau^{t,x})\geq\phi(\tau,X_\tau^{t,x})$ and for each $r\in[t,\tau],$
$$\widetilde{G}(r,\phi(r,X_r^{t,x}),\sigma^*\nabla_x\phi(r,X_r^{t,x}))\geq-\frac{\partial\phi}{\partial t}(r,X_r^{t,x})-{\cal{L}}\phi(r,X_r^{t,x})+C.\eqno(4.7)$$
Then, since $({\tilde{Y}}_s,{\tilde{Z}}_s)\in{\mathcal{S}}_D^\infty
\times {\cal{H}}^\infty,$ by (4.6), (4.7) and Proposition 3.11(ii) (strict comparison), we get $\hat{Y}_t>\tilde{Y}_t$. Since ${\bar{Y}}_t\geq\hat{Y}_t$, it follows that $v(t,x)=\bar{Y}_t>\tilde{Y}_t=\phi(t,x),$ which contradicts the condition that $v(t,x)=\phi(t,x).$ Thus, $v(t,x)$ is a viscosity supersolution of (4.1).

Step 2.  We will show that $v(t,x)$ is a viscosity subsolution of (4.1). This proof is similar to Step 1. Suppose that $v(t,x)$ is not a viscosity subsolution of (4.1). This means that there exists $(t,x,\phi)\in[0,T]\times \textbf{R}^d\times C^{1,2}_D([0,T]\times \textbf{R}^d)$ satisfying the condition that $\phi(t,x)=v(t,x)$ and that $(t,x)$ is a local minimum point of $\phi-v$, such that $v(t,x)>h(t,x)$ and
$$\frac{\partial\phi}{\partial t}(t,x)+{\cal{L}}\phi(t,x)+\widetilde{G}(t,\phi(t,x),\sigma^*\nabla_x\phi(t,x))<0.$$
By continuity, there exist $\beta\in(0,T-t],$ $c>0,$ and $C>0$ such that for each $s\in[t,t+\beta]$ and $y\in[x-c,x+c]$, we have $v(s,y)\leq\phi(s,y),$ $$v(s,y)\geq h(s,y)+C,\ \ \textmd{and}\ \
\frac{\partial\phi}{\partial t}(s,y)+{\cal{L}}\phi(s,y)+\widetilde{G}(s,\phi(s,y),\sigma^*\nabla_x\phi(s,y))\leq -C.\eqno(4.8)$$
We define a stopping time $\tau=\inf\{s\geq t; |X_s^{t,x}-x|\geq c\}\wedge(t+\beta),$ then $t<\tau \leq t+\beta$, and the ranges of $v({s\wedge\tau},X_{s\wedge\tau}^{t,x})_{s\in[t,t+\beta]}$, $h({s\wedge\tau},X_{s\wedge\tau}^{t,x})_{s\in[t,t+\beta]}$ and $\phi({s\wedge\tau},X_{s\wedge\tau}^{t,x})_{s\in[t,t+\beta]}$ are all included in a closed subset of $D$. By (4.3) and (4.8), we have, for each $s\in[t,\tau],$
$$Y^{t,x}_s=v(s,X^{t,x}_s)\geq h(s,X^{t,x}_s)+C,$$
which implies that for each $s\in[t,\tau],$ $dK^{t,x}_s=0.$ Then by (4.2) and (4.3), we get that the following BSDE:
$$\bar{Y}_s=v(\tau,X_\tau^{t,x})+\int_s^{t+\beta} 1_{\{r\leq\tau\}}\widetilde{G}(r,\bar{Y}_r,\bar{Z}_r)dr-\int_s^{t+\beta}\bar{Z}_rdB_r,\ \ s\in[t,t+\beta],\eqno(4.9)$$
has a solution $(\bar{Y}_s,\bar{Z}_s)=({Y}_{s\wedge\tau}^{t,x},1_{\{s\leq\tau\}}{Z}_s^{t,x})\in{\mathcal{S}}_D^\infty
\times{{H}}^2.$ Applying It\^{o}'s formula to $\phi(s,X_s^{t,x})$ for $s\in[t,\tau]$, we get that the following BSDE:
$${\tilde{Y}}_s=\phi(\tau,X_\tau^{t,x})+\int_s^{t+\beta} -1_{\{r\leq\tau\}}\left(\frac{\partial\phi}{\partial t}(r,X_r^{t,x})+{\cal{L}}\phi(r,X_r^{t,x})\right)dr-\int_s^{t+\beta}{\tilde{Z}}_rdB_r,\ \ s\in[t,t+\beta],\eqno(4.10)$$
has a solution $({\tilde{Y}}_s,{\tilde{Z}}_s)=(\phi({s\wedge\tau},X_{s\wedge\tau}^{t,x}), 1_{\{s\leq\tau\}}\sigma^*\nabla_x\phi(s,X_s^{t,x}))\in{\mathcal{S}}_D^\infty
\times {\cal{H}}^\infty.$ By (4.8), we have $v(\tau,X_\tau^{t,x})\leq\phi(\tau,X_\tau^{t,x})$ and for each $r\in[t,\tau],$
\begin{eqnarray*}
-\frac{\partial\phi}{\partial t}(r,X_r^{t,x})-{\cal{L}}\phi(r,X_r^{t,x})\geq\widetilde{G}(r,v(r,X_r^{t,x}),\sigma^*\nabla_x\phi(r,X_r^{t,x})+C.
\end{eqnarray*}
Then, since $({\tilde{Y}}_s,{\tilde{Z}}_s)\in{\mathcal{S}}_D^\infty
\times {\cal{H}}^\infty,$ it follows from that (4.9), (4.10) and Proposition 3.11(i) that $v(t,x)=\bar{Y}_t<\tilde{Y}_t=\phi(t,x),$ which contradicts the condition that $v(t,x)=\phi(t,x).$ Thus, $v(t,x)$ is a viscosity subsolution of (4.1).

The proof is complete.\ \ $\Box$
\subsection{An optimal stopping problem for the payoff of American options}
From the definition of $u_f$, we observe that $u_f$ is a utility function (concave and strictly increasing), when $f$ is nonpositive. More specially, from Example 2.1, we get that
\begin{flushleft}
\ \ \ \ \ \ \ \ $u_f$ is an exponential utility, when $f(y)=-\beta,\ \beta>0;$\\
\ \ \ \ \ \ \ \ $u_f$ is a power utility, when $f(y)=-\frac{\beta}{y},\ \beta>0,\ \beta\neq\frac{1}{2};$\\
\ \ \ \ \ \ \ \ $u_f$ is a logarithmic utility, when $f(y)=-\frac{1}{2y}.$
\end{flushleft}
From Lemma A.1(iii), we observe that $-2f(x)$ is the Arrow-Pratt risk aversion coefficient of the utility function $u_f(x)$ (see \cite[Definition 2.45 and Example 2.46]{FS}).

The holder of an American option has the right to exercise the option at any stopping time $\tau\in{\cal{T}}_{0,T}$ according to the expected utility of its payoff. Let $f$ be nonpositive and let the payoff of the American option be described by $\eta_t:=L_t1_{\{t<T\}}+\xi1_{\{t=T\}}$ with $u_f(\eta_t)\in {\cal{S}}^p$. For $\tau\in{\cal{T}}_{0,T}$, the conditional expected utility of $\eta_{\tau}$ is given by $E[u_f(\eta_\tau)|{\cal{F}}_t]$. Then, by Corollary 3.5 and Remark 3.7, the RBSDE$(f(y)|z|^2,\xi,L_t)$ has a unique solution $(y_t,z_t,k_t)$ such that $u_f(y_t)$ is the maximal conditional expected utility of the payoff $\eta_t$, i.e.,
$$u_f(y_t)=\textrm{ess}\sup_{\tau\in{\cal{T}}_{t,T}} E[u_f(\eta_\tau)|{\cal{F}}_t],\ \ t\in[0,T].\eqno(4.11)$$
For $t\in[0,T]$, set $\sigma_t^\ast:=\inf\{s\geq t,y_s=\eta_s\}\wedge T.$ Since $dk_s=0$ on $[t,\sigma_t^\ast]$, by Lemma A.1(ii) and (4.11), we get that $\sigma_t^\ast$ is an optimal stopping time such that
$$E[u_f(\eta_{\sigma_t^\ast})|{\cal{F}}_t]=\textrm{ess}\sup_{\tau\in{\cal{T}}_{s,T}} E[u_f(\eta_\tau)|{\cal{F}}_t].$$

More generally, we consider the notion of $F$-evaluation ${\cal{E}}^F_{s,t}[\cdot]$, which is a nonlinear evaluation introduced by \cite[Definition 3.1]{Peng}. For $\sigma\in{\cal{T}}_{0,T}$, let $(\hat{y}_t,\hat{z}_t)$ be the solution to the BSDE$(F,u_f(\eta_\sigma))$ with terminal time $\sigma$ such that $\hat{y}_t\in{\cal{S}}^p$. Then, we denote the $F$-evaluation of the utility of $\eta_\sigma$ at time $t$ by $${\cal{E}}^{F}_{t,\sigma}[u_f(\eta_\sigma)]:=\hat{y}_t,\ \  t\in[0,T].$$
By \cite[Theorem 5.9]{EQ}, the RBSDE$(F,u_f(\xi),u_f(L_t))$ has a unique solution $(\tilde{y}_t,\tilde{z}_t,\tilde{k}_t)$ such that
$$\tilde{y}_t=\textrm{ess}\sup_{\tau\in{\cal{T}}_{t,\tau}}{\cal{E}}^{F}_{t,\tau}[u_f(\eta_\tau)],\ \ t\in[0,T].$$
Let us further assume that $e^{\gamma t}(u_f(\eta_t)\vee0)+\delta t)\in{\mathcal{C}}_{V}$ with $u_f(\eta_t)\in{\mathcal{S}}^p$. By Theorem 3.3 and (3.5), the RBSDE$(G^F_f+f(y)|z|^2,\xi,L_t)$ has a unique solution $(y_t,z_t,k_t)$ such that $$u_f(y_t)=\tilde{y}_t=\textrm{ess}\sup_{\tau\in{\cal{T}}_{t,\tau}} {\cal{E}}^{F}_{t,T}[u_f(\eta_\tau)],\ \ t\in[0,T].\eqno(4.12)$$
For $t\in[0,T]$, set $\sigma_t^\ast:=\inf\{s\geq t, y_s=\eta_s\}\wedge T$. Since $dk_s=0$ on $[t,\sigma_t^\ast]$, by Lemma A.1(ii) and (4.12), we get that $\sigma_t^\ast$ is an optimal stopping time such that
$${\cal{E}}^{F}_{t,\sigma_t^\ast}[u_f(\eta_{\sigma_t^\ast})]=\textrm{ess}\sup_{\tau\in{\cal{T}}_{t,\tau}}{\cal{E}}^{F}_{t,\tau}[u_f(\eta_\tau)].$$
\\\\

\begin{center}
\textbf{Appendix: A}
\end{center}

Let $f\in L_{1,loc}(D)$ be given. Given $\alpha\in D,$ we define the following:
$$u_f(x):=\int^x_\alpha\exp\left(2\int^y_\alpha f(z)dz\right)dy,\ \ x\in D.$$
The following properties for $u_f(x)$ come from \cite[Lemma A.1]{B17} and \cite[Lemma 2.1]{Zheng1}.\\\\
\textbf{Lemma A.1}\ \emph{The following properties of $u_f(x)$ hold:}

\emph{(i) $u_f(x)\in W^2_{1,loc}(D)$, in particular, $u_f(x)\in C^1(D)$;}

\emph{(ii) $u_f(x)$ is strictly increasing;}

\emph{(iii) $u_f''(x)-2f(x)u_f'(x)=0,\ a.e.$ on $D$;}

\emph{(iv) $u_f^{-1}(x)\in W^2_{1,loc}(V)$, in particular, $u_f^{-1}(x)\in C^1(V)$ and is strictly increasing;}

\emph{(v) If $l\in L_{1,loc}(D)$ and $l(x)\leq f(x), a.e.$, then for every $x\in D,$ $u_l(x)\leq u_f(x)$.}\\

To deal with discontinuous generators, an It\^{o}-Krylov formula was established in \cite[Theorem 2.1]{B17}. To treat our situation conveniently, we give in Lemma A.2 this It\^{o}-Krylov formula for semimartingales with values in $D$, using a similar argument.\\\\
\textbf{Lemma A.2}\ \ \emph{Let $Y_t=Y_T+\int_t^Th(s)ds+A_T-A_t-\int_t^TZ_sdB_s, t\in[0,T],$
where $Y_t\in{\cal{C}}_D$, $A_t\in {\cal{A}},$ $Z_t\in H^2$ and $h(t)$ is a progressively measurable process satisfying $\int_0^T|h(t)|dt<\infty.$ Then the following two statements hold:}

\emph{(i) (Krylov's estimate) Let $\{O_n\}_{n\geq1}$ be a sequence of closed intervals such that for each $n\geq1,$ $O_n\subset O_{n+1}$, and $Y_0\in O_1,$ $\cup_{n\geq1}O_n=D.$ Let $\tau_n^1=\inf\{t\geq0,Y_t\notin O_n\}\wedge T$, $\tau_n^2=\inf\{t\geq0,(A_t+\int_0^t|h(s)|ds+\int_0^t|Z_s|^2ds)>n\}\wedge T$ and $\tau_n=\tau_n^1\wedge\tau_n^2.$ Then for each nonnegative $\psi\in L_{1,loc}(D)$, we have } $$E\left[\int_0^{T\wedge\tau_n}\psi(Y_t)|Z_t|^2dt\right]\leq (2n+2\lambda(O_n))\int_{O_n}\psi(x)dx.$$

\emph{(ii) (It\^{o}-Krylov formula) For each $u\in W^2_{1,loc}(D)$, we have}
 $$u(Y_t)=u(Y_0)+\int_0^tu'(Y_s)dY_s+\frac{1}{2}\int_0^tu''(Y_s)|Z_s|^2ds, \ \ t\in[0,T].$$
\emph{Proof.} It is clear that there exists a sequence of closed intervals $\{O_n\}_{n\geq1}$ such that for any $n\geq1,$ $O_n\subset O_{n+1}$, and $Y_0\in O_1,$ $\cup_{n\geq1}O_n=D.$ Let $\tau_n^1=\inf\{t\geq0,Y_t\notin O_n\}\wedge T$, $\tau_n^2=\inf\{t\geq0,(A_t+\int_0^t|h(s)|ds+\int_0^t|Z_s|^2ds)>n\}\wedge T$ and $\tau_n=\tau_n^1\wedge\tau_n^2.$

Proof (i): Using Tanaka's formula, for $a\in D$, we have
\begin{align*}
(Y_{t\wedge\tau_n}-a)^-=(Y_0-a)^-&+\int_0^{t\wedge\tau_n}1_{\{Y_s\leq a\}}h(s)ds+\int_0^{t\wedge\tau_n}1_{\{Y_s\leq a\}}dA_s\\
&-\int_0^{t\wedge\tau_n}1_{\{Y_s\leq a\}}Z_sdB_s+\frac{1}{2}L_{t\wedge\tau_n}^a(Y),\ \ t\in[0,T].
\end{align*}
Since
$$|(Y_{t\wedge\tau_n}-a)^--(Y_0-a)^-|\leq|Y_{t\wedge\tau_n}-Y_0|\leq\lambda(O_n),$$
we deduce that for each $a\in D$ and $t\in[0,T],$
$$E(L_{t\wedge\tau_n}^a(Y))\leq2n+2\lambda(O_n).$$
It follows from occupation times formula that for each nonnegative $\psi\in L_{1,loc}(D)$,
\begin{eqnarray*}
E\left[\int_0^{T\wedge\tau_n}\psi(Y_t)|Z_t|^2dt\right]&=&E\left[\int_0^{T\wedge\tau_n}\psi(Y_t)d\langle Y,Y\rangle_t\right]\\
&=&\int_{O_n}\psi(a)E(L_{T\wedge\tau_n}^a(Y)da\\
&\leq&(2n+2\lambda(O_n))\int_{O_n}\psi(a)da
\end{eqnarray*}

Proof of (ii): For any $u\in W^2_{1,loc}(D)$, using a convolution method, we can find a sequence $\{u_m\}_{m\geq1}$ in $C^2(D)$ such that

 1) $u_m$ converges uniformly to $u$ in $O_n$;

 2) $u'_m$ converges uniformly to $u'$ in $O_n$;

 3) $u''_m$ converges to $u'$ in $L^1(O_n)$.

Using It\^{o}'s formula, we have
 $$u_m(Y_{t\wedge\tau_n})=u_m(Y_0)+\int_0^{t\wedge\tau_n}u_m'(Y_s)dY_s+\frac{1}{2}\int_0^{t\wedge\tau_n}u_m''(Y_s)|Z_s|^2ds, \ \ t\in[0,T].$$
Then, by (i) and the same arguments as in the proof of \cite[Theorem 2.1]{B17}, we get (ii). \ \ $\Box$

\end{document}